\documentclass[oneside, draft, 11pt]{amsart}  
\newtheorem{theorem}{Theorem}[section]
\newtheorem{prop}[theorem]{Proposition}
\newtheorem{lemma}[theorem]{Lemma}
\theoremstyle{definition}

\theoremstyle{remark}
\newtheorem{example}[theorem]{Example}
\newcommand{\mbf}{\mathbf}
\newcommand{\mcl}{\mathcal}
\newcommand{\ZZ}{\mathbb{Z}}
\newcommand{\datum}{6/20/'17} 
\newcommand{\isomto}{\cong}
\numberwithin{equation}{section}
\usepackage{amssymb, latexsym, amscd}
\usepackage{epsfig}
\usepackage{MnSymbol}

\begin{document}

\title[S.P.\ Ellis, Joining and Independence, \datum]{Joining and Independence in Concurrence Topology}
       
\author[S.~P.~Ellis]{Steven P.\ Ellis}
\date{\datum}

\address{Unit 42, NYSPI \\
1051 Riverside Dr. \\
New York, NY 10032 \\
U.S.A.}

\email{ellisst@nyspi.columbia.edu}


        
\thanks{This research is supported in part by United States PHS grants MH46745, MH60995, and MH62185.}

\thanks{\datum}

 \begin{abstract}
``Concurrence topology'' (Ellis and Klein \emph{Homology, Homotopy, and Applications,} \textbf{16}) is a TDA method for binary data. The idea is to construct a filtration consisting of Dowker complexes then compute persistent homology. Persistent classes correspond to a form of negative statistical association among the variables.

Suppose we have two groups of binary variables each displaying negative association, manifested in nontrivial concurrence homology in dimensions $p$ and in one group and $q$ in the other \emph{when the groups of variables are considered individually.} Suppose, however, that the two \emph{groups} of variables are statistically independent of each other. Now combine the two groups of variables and suppose the sample size is large. Then representative cycles, one from each group of variables, will combine to produce a cycle in dimension $p+q+1$. This is a chain level phenomenon, but we show it has a signature in homology. Looking for this signature can be used to study the dependence among groups of variables. 
 \end{abstract}
 
 \maketitle

\section{Introduction}   \label{S:intro}
Let $K$ and $L$ be finite simplicial complexes. Let $\sigma \in K$ with vertices $v_{1}, \ldots, v_{m}$. Let $\tau \in K$ with vertices $v_{m+1}, \ldots, v_{n}$. \emph{We assume} that for every $\sigma$ and $\tau$ the sequence 
\newline
$v_{1}, \ldots, v_{m}, v_{m+1}, \ldots, v_{n}$ is geometrically independent (Munkres \cite[p.\ 12]{jrM84}). Let $\sigma \ast \tau$ be the simplex with vertices $v_{1}, \ldots, v_{m}, v_{m+1}, \ldots, v_{n}$  (Munkres \cite[p.\ 368]{jrM84}). Assume that the join, $K \ast L$, exists. Recall that in that case $K \ast L$ is the complex consisting 
of all $\sigma \ast \tau$ with $\sigma \in K$ and $\tau \in L$ and \emph{all faces of such.} Thus, $\sigma, \tau \in K \ast L$. 
 
Assume 
	\begin{equation} \label{E:K.L.disjoint}
		K \cap L = \varnothing.
	\end{equation}

Consider the space $X$ that is a union 
	\begin{equation}  \label{E:X.is.union.of.prods}
		X := \bigcup_{i=1}^{N} \sigma_{i} \times \tau_{i},
	\end{equation}
where $\sigma_{i} \in K$, $\tau_{i} \in L$ are simplices of $K$ and $L$, resp. The $\sigma_{i}$'s do not have to be distinct. Nor do the $\tau_{i}$'s. However, we assume
	\begin{equation} \label{E:sigma.tau.pairs.distinct}
		i,j = 1, \ldots, N \text{ and } i \neq j \text{ imply } (\sigma_{i}, \tau_{i}) \neq (\sigma_{j}, \tau_{j}).
	\end{equation}
We assume further that every facet (i.e., maximal simplex) of $K$ appears at least once as a $\sigma_{i}$. Ditto for $L$. Thus, there is a map $g : \{1, \ldots, N \} \to K \times L$ that associates to every index $i = 1, \ldots, N$ a pair $g(i) = ( \sigma_{i}, \tau_{i} ) \in K \times L$. The projection onto the first (second) factor in $g \bigl( \{ 1, \ldots, N \} \bigr)$ includes the set of facets of $K$ ($L$, resp.). The set of vertices $K^{(0)}$ and $L^{(0)}$ are disjoint. In particular, 
	\begin{equation} \label{E:sigma.i.tau.i.disjoint}
		\sigma_{i} \cap \tau_{i} = \varnothing \text{ for every } i = 1, \ldots, N.
	\end{equation}

In this paper except where noted \emph{all homology is with $\ZZ/2$ coefficients.} Let $\mcl{C}(K) = \bigl\{ \ldots \xrightarrow{\partial} C_{2}(K) \xrightarrow{\partial} C_{1}(K) \xrightarrow{\partial} C_{-1}(K) = \ZZ \bigr\}$ and $\mcl{C}(L) = \bigl\{ \ldots \xrightarrow{\partial} C_{2}(L) \xrightarrow{\partial} C_{1}(L) \xrightarrow{\partial} C_{-1}(L) = \ZZ \bigr\}$ be the \emph{augmented} integer simplicial chain complexes of $K$ and $L$, resp.\. 
If $c = \sigma_{1} + \cdots + \sigma_{m}$ is a $p$-chain in $\mcl{C}(K) \otimes \ZZ/2$ and $d = \tau_{1} + \cdots + \tau_{n}$ is a $q$-chain in $\mcl{C}(L) \otimes \ZZ/2$ ($p, q \geq 0$), define the chain $c \ast d$ in $K \ast L$ to be the $(p+q+1)$-chain 
$\sum_{i=1}^{m} \sum_{j=1}^{n} \sigma_{i} \ast \tau_{j}$. It is easy to see that if $c$ and $d$ are cycles then so is $c \ast d$. Denote the collection of all such $c \ast d$ together with $\mcl{C}(K)$ and $\mcl{C}(L)$ by $\mcl{C}(K) \ast \mcl{C}(L)$. 
 
Write $\mcl{C}(K) \otimes_{\ZZ/2} \mcl{C}(L) := \bigl[ \mcl{C}(K) \otimes \ZZ/2 \bigr] \otimes_{\ZZ/2} \bigl[ \mcl{C}(L) \otimes \ZZ/2 \bigr]$. As usual, we define $\dim (\sigma \otimes \tau) = \dim \sigma + \dim \tau$. We think of the ``chain'' 1 as belonging to $C_{-1}(K)$ or $C_{-1}(L)$ so $\dim (\sigma \otimes 1) = \dim \sigma -1$ and $\dim (1 \otimes \tau) = \dim \tau -1$. Thus, $\mcl{C}(K) \otimes_{\ZZ/2} \mcl{C}(L)$ has non-trivial homology in dimension -2: 
$\bigl[ \mcl{C}(K) \otimes_{\ZZ/2} \mcl{C}(L) \bigr]_{-2} = \{ 1 \otimes 1 \}$. For $\sigma \in K$, $\tau \in L$, 
define $\partial (\sigma \otimes \tau) := (\partial \sigma) \otimes \tau + \sigma \otimes (\partial \tau)$.

The proof of the following is easy.

  \begin{lemma}
For $r \geq 0$ define $f : C_{r+1}(K \ast L) \otimes \ZZ/2 \to \bigl[ \mcl{C}(K) \otimes_{\ZZ/2} \mcl{C}(L) \bigr]_{r}$ 
by $f(\sigma \ast \tau) := \sigma \otimes \tau$ ($\sigma \in K$, $\tau \in L$). If $r \geq -1$ and $\sigma \in K$, $\tau \in L$ with $\dim \sigma = \dim \tau = r+1$, define $f(\sigma) := \sigma \otimes 1 \in C_{r+1}(K) \otimes_{\ZZ/2} C_{-1}(L)$
and $f(\tau) = 1 \otimes \tau \in C_{-1}(K) \otimes_{\ZZ/2} C_{r+1}(L)$ 
Finally, with $r = -2$, set $f(1) := 1 \otimes 1 \in C_{-1}(K) \otimes_{\ZZ/2} C_{-1}(L)$. 
Then $f$ is a chain map.
  \end{lemma}

Note that the lemma does not seem to hold with $\ZZ/2$ replaced by an arbitrary field: There can be a problem with signs. 

Now let $\alpha \in \tilde{H}_{p}(K)$ and $\beta \in \tilde{H}_{q}(L)$ be non-trivial ($p, q \geq 0$). Let $w$ and $z$ be cycles representing $\alpha$ and $\beta$, resp.\. Then $w \ast z \in C_{p+q+1}(K \ast L) \otimes \ZZ/2$. In fact, $w \ast z$ is obviously a cycle (with $\ZZ/2$ coefficients). Call $w \ast z$ a ``join cycle''. We are interested in such cycles for reasons explained in section \ref{S:stat.interp}. Thus, our interest is at the chain, not homology level. We wish to be able to see at the homology level, signs of the presence of join cycles.

Suppose $w \ast z$ in fact lies in the subgroup 
$C_{p+q+1}(M) \otimes \ZZ/2$. Let $i : \mcl{C}(M) \otimes \ZZ/2 \to \mcl{C}(K \ast L) \otimes \ZZ/2$ 
be inclusion\footnote{More precisely, let $p \geq 0$ and let $i' : C_{p}(M) \to C_{p}(K \ast L)$ be inclusion. Let $G_{2}$ be the subgroup of $C_{p}(K \ast L)$ spanned by the $p$-simplices in $C_{p}(K \ast L)$ \emph{not} in $M$. 
Then $C_{p}(K \ast L) = C_{p}(M) \oplus G_{2}$. Therefore, by Munkres \cite[Top of p.\ 24]{jrM84}, we have
$C_{p}(K \ast L)/C_{p}(M) = G_{2}$ so is free abelian. Therefore, by Munkres \cite[Corollary 23.2, p.\ 132]{jrM84}, the exact sequence
$0 \to C_{p}(M) \xrightarrow{i'} C_{p}(K \ast L) \to C_{p}(K \ast L)/C_{p}(M) \to 0$ splits. 
Hence, by Munkres \cite[Theorem 50.4, p.\ 302]{jrM84}, the homomorphism $i := i' \otimes \ZZ/2$ is injective.} 
 and suppose $i(w \ast z)$ represents a non-trivial class in $\tilde{H}_{p+q+1}(M)$. Consider the composition of chain maps
  \begin{equation*}
	\mcl{C}(M) \otimes \ZZ/2 \xrightarrow{i} \mcl{C}(K \ast L) \otimes \ZZ/2 \xrightarrow{f} \mcl{C}(K) \otimes_{\ZZ/2} \mcl{C}(L).
  \end{equation*}
Now apply the homology functor. Then, by theorem \ref{E:Kunneth.with.field}, we have this composition of homomorphisms:
  \begin{equation*}
	\tilde{H}_{\ast}(M) \xrightarrow{i_{\ast}} \tilde{H}_{\ast}(K \ast L) 
	  \xrightarrow{f_{\ast}} \tilde{H}_{\ast} \bigl[ \mcl{C}(K) \otimes_{\ZZ/2} \mcl{C}(L) \bigr] 
	    \xrightarrow{Kun^{-1}} \tilde{H}_{\ast}(K) \otimes_{\ZZ/2} \tilde{H}_{\ast}(L).
  \end{equation*}

This composition of homomorphisms takes $\{ w \ast z \} \in \tilde{H}_{\ast}(M)$ to $\{ w \} \otimes \{ z \} = \alpha \otimes \beta$, which is nontrivial by assumption. Therefore, $i_{\ast} \bigl[ \{ w \ast z \} \bigr] \in \tilde{H}_{\ast}(K \ast L)$ is non-trivial. To sum up:


  \begin{theorem}  \label{T:i.w*z.nontriv}
Let $K$ and $L$ be simplicial complexes, perhaps abstract.  
Let $\alpha \in \tilde{H}_{p}(K)$ and $\beta \in \tilde{H}_{q}(L)$ be non-trivial ($p, q \geq 0$). Let $w$ and $z$ be cycles representing $\alpha$ and $\beta$, resp. Suppose $w \ast z \in C_{p+q+1}(M) \otimes \ZZ/2$ and $w \ast z$ represents a non-trivial class in $\tilde{H}_{p+q+1}(M)$. 
Let $i : \mcl{C}(M) \otimes \ZZ/2 \to \mcl{C}(K \ast L) \otimes \ZZ/2$ be inclusion. Then $i(w \ast z)$ represents a non-trivial class 
in $\tilde{H}_{p+q+1}(K \ast L)$. 
  \end{theorem}
  
Here is a converse:
  \begin{prop}
Let $K$ and $L$ be simplicial complexes, perhaps abstract. Let $w$ and $z$ be cycles in $\mcl{C}(K) \otimes \ZZ/2$ 
and $\mcl{C}(L) \otimes \ZZ/2$, resp.\ s.t.\ $w \ast z$ represents non-trivial reduced homology in $K \ast L$. Then $w$ and $z$ represent non-trivial reduced homology in $K$ and $L$, resp.
  \end{prop}
  \begin{proof}
Suppose $i(w \ast z)$ represents non-trivial reduced homology in $K \ast L$ but $w$, say, does not represent non-trivial reduced homology in $K$. Then there exists a chain $c \in \mcl{C}(K)$ s.t.\ $\partial c = w$. Write $w = \sum_{i=1}^{m} \sigma_{i}$, $z = \sum_{j=1}^{n} \tau_{j}$, and $c = \sum_{k=1}^{\ell} \zeta_{k}$. Let $d = c \ast z$, a chain in $K \ast L$. Then
    \begin{align*}
        \partial d &= \sum_{j,k} ( \partial \zeta_{k} \ast \tau_{j} + \zeta_{k} \ast \partial \tau_{j} ) \\
          &=  \sum_{j} \left( \sum_{k} \partial \zeta_{k} \right) \ast \tau_{j} 
            + \sum_{k} \zeta_{k} \ast \left( \sum_{j=1}^{n} \partial \tau_{j} \right) \\
          &=  \sum_{j} (\partial c) \ast \tau_{j} + \sum_{k} \zeta_{k} \ast ( \partial z ) \\
          &=  \sum_{j} w \ast \tau_{j} + \sum_{k} \zeta_{k} \ast 0 \\
          &=  w \ast z.
    \end{align*}
This contradicts the assumption that $w \ast z$ represents non-trivial homology, i.e., it does not bound.
  \end{proof}

  \begin{example}
  1) Let $K := \{ a, b \}$,  $L := \{A, B, AB \}$, and $M := \{ a, b, A, B, AB, aA, aB, bA, bB \}$. Then $w := a + b$ is a 0-cycle in $K$ that represents non-trivial (reduced) homology in dimension 0. On the other hand, $z := A + B$ is a 0-cycle in $L$ that \emph{does not} represent non-trivial homology in dimension 0: $z = \partial AB$. Still, $y := w \ast z \in Z_{1}(M)$ represents non-trivial homology in dimension 1. But $i(y)$ does not represent non-trivial homology in $K \ast L$: $i(y) = \partial(aAB + bAB)$. Thus, that $w \ast z$ represents a non-trivial class in $\tilde{H}_{p+q+1}(M)$ does not automatically imply that $i_{\ast} \bigl[ \{ w \ast z \} \bigr]$ is non-trivial.
  
  2) Let $K$ consist of three vertices: $K := \{ a, b, c \}$. Similarly, let $L$ consist of three vertices: $L := \{ A, B, C \}$. Let $M$ be the sub complex of $K \ast L$ whose facets are $aA$, $Ab$, $bB$, $Bc$, $cC$, and $Ca$. The cycle $z := aA + Ab + bB + Bc + cC + Ca$ in $C_{1}(M)$ represents non-trivial homology in dimension 1. Since $K \ast L$ is one-dimensional, $z$ also represents non-trivial homology in $K \ast L$. However, $z$ is not the ``join'' of two cycles, one from $K$ and one from $L$. Obviously, this idea generalizes so that if $n \geq 2$ and $K$ and $L$ each consist of $n$ vertices then one can get a 1-cycle $z$ with $2n$ terms. 

 3) Consider the complex $K$ with facets $a, bd, cd, bc$ and the complex $L$ with facets $D$, $AB$, $BC$, $AC$. $K$ and $L$ each have nontrivial reduced homology in dimensions 0 and 1. Combine them into a complex $M$ with facets 
$aAB, aAC, aBC, AbB, bBc, BcC, cCd, AbC, bCd, bdD, cdD, bcD$.

The sum of the facets of $M$, \emph{viz.} $z := aAB + aAC + aBC + AbB + bBc + BcC + cCd + AbC + bCd + bdD + cdD + bcD$ is a two cycle that represents a non-trivial homology class, $\alpha$. Since $\dim M = 2$, $z$ is the only representative of $\alpha$ and $\alpha$ is the only class in $H_{2}(M)$. Let $i : M \hookrightarrow K \ast L$ be inclusion. It turns out that $i_{\ast}(\alpha) \neq 0$ in $H_{2}(K \ast L)$.

Now, $z$ is not the join of cycles from $K$ and $L$. However, $z = c_{1} + c_{2}$, where
	\begin{equation}
		c_{1} =  aAB + aAC + aBC + AbB + BcC + AbC \text{ and } c_{2} = bBc + cCd + bCd + bdD + cdD + bcD.
	\end{equation}
$c_{1}$ and $c_{2}$ are \emph{almost}, but not quite, the joins of cycles, one from $K$ and one from $L$. 
In fact, neither $c_{1}$ nor $c_{2}$ is even a cycle.
  \end{example}

\section{Statistical interpretation} \label{S:stat.interp}
Now suppose $V_{1}, \ldots, V_{n}$ are dichotomous (binary) variables and it is natural to assign them to two groups, $V_{1}, \ldots, V_{m}$ and $V_{m+1}, \ldots, V_{n}$. (Generalizing to more than two groups is straight forward.) Suppose we collect $T$ observations $X_{1}, \ldots, X_{T}$ (collectively, $\mbf{X}$), on these variables and apply the ``concurrence topology'' (Ellis and  Klein \cite{spEaK14.ConcurTopolfMRI}) method to those data. The observations, $X_{1}, \ldots, X_{T}$, do not have to be independent, but suppose that the multiple time series 
$\mbf{X}$ ``mixes'' sufficiently that the observations exhibit independent-like statistics 
(Brillinger \cite[Section 1.3]{drB01.BrillingersTimeSeriesBook}). Recall that in that method one constructs an abstract descending filtered simplicial complex $\mcl{M}$. (Each frame in the filtration is a ``Dowker complex'' (Dowker \cite{Dch52.DowkerComplex}.) The vertices are just $V_{1}, \ldots, V_{n}$. The vertices $V_{i_{1}} \ldots V_{i_{k}}$ define a simplex 
of $\mcl{M}$ if and only if in some observation we have $V_{i_{1}} = 1, \ldots, V_{i_{k}} = 1$. Denote the subcomplexes that only use $V_{1}, \ldots, V_{m}$ and $V_{m+1}, \ldots, V_{n}$ by $\mcl{K}$ and $\mcl{L}$, resp. 

Suppose the random vectors $(V_{1}, \ldots, V_{m})$ and $(V_{m+1}, \ldots, V_{n})$ are independent of each other. Because of the independence, for large $T$ if $\sigma \in \mcl{K}$ and $\tau \in \mcl{L}$ are simplices that each occur with positive probability in $\mcl{K}$ and $\mcl{L}$, resp.\ then the simplex $\sigma \ast \tau$ will eventually appear in $\mcl{M}$. Thus, the complex $\mcl{M}$ will sort of look like $\mcl{K} \ast \mcl{L}$.

But suppose that, while the random vectors $(V_{1}, \ldots, V_{m})$ and $(V_{m+1}, \ldots, V_{n})$ are independent of each other, the coordinates of the first vector, $V_{1}, \ldots, V_{m}$, are \emph{not} independent of each other. Specifically, suppose that the joint distribution of $V_{1}, \ldots, V_{m}$ is such that with ``high probability'' the complex $\mcl{K}$ has a persistent class $\beta$ in dimension $p \geq 0$, say. Similarly, suppose that with high probability $\mcl{L}$ has a persistent class $\gamma$ in dimension $q \geq 0$. Let $w := \sum_{i=1}^{s} \sigma_{i}$ be a representative cycle for $\beta$ and $z := \sum_{i=1}^{t} \tau_{i}$ be a representative cycle for $\gamma$.   

Then, because of independence between $V_{1}, \ldots, V_{m}$ and $V_{m+1}, \ldots, V_{n}$ and the ``mixing'' behavior of $X_{1}, \ldots, X_{T}$, there will appear in $\mbf{X}$ every combination $\sigma_{i} \ast \tau_{j}$.    
Specially if ``with high probability'' the classes $\beta$ and $\gamma$ have long lifespans, there will eventually be frequency levels at which one observes the join $w \ast z$. 

Hence, by theorem \ref{T:i.w*z.nontriv} in some frames the homomorphism 
$i_{\ast} : C_{p+q+1}(M) \otimes \ZZ/2 \to C_{p+q+1}(K \ast L) \otimes \ZZ/2$ will be nontrivial. This we will take as a sign of possible independence.

In fact, this line of reasoning inspires a weakening of the notion of join cycle. Suppose $w \ast z$ is a join cycle. Suppose $\xi \in K \ast L$ whose vertices belong to $(V_{1}, \ldots, V_{m}, V_{m+1}, \ldots, V_{n})$ and $\sigma_{i}$ and $\tau_{j}$ are faces of $\xi$. Suppose $c$ is a $p+q+1$ chain whose terms are faces of $\xi$ and, roughly, speaking, $\partial \xi + c = \sigma_{i} \ast \tau_{j}$. Then replacing $\sigma_{i} \ast \tau_{j}$ in $w \ast z$ by $c$ we get a cycle, $x$, that we will regard as a kind of generalized join cycle. The argument given above implies that ``in the fullness of time'' we should frequently get the cycle $x$.

In this paper, we use this line of reasoning to inspire a method for assessing the weakness of association between two sets of variables, $V_{1}, \ldots, V_{m}$ and $V_{m+1}, \ldots, V_{n}$. The more classes $\gamma \in H_{\ast}(M)$ s.t.\ 
$\partial_{\ast} \circ i_{\ast}(\gamma)$
is nontrivial and the longer its lifespans, the more ``independence-like'' behavior we ascribe to the pair $V_{1}, \ldots, V_{m}$ and $V_{m+1}, \ldots, V_{n}$. This is reminiscent of ``canonical correlation analysis'' in classical multivariate analysis 
(Johnson and Wichern \cite[Chapter 10]{raJdwW92}). However, we only examined the problem of studying the dependence among \emph{two} groups of variables for convenience. The same ideas easily apply to any number of groups of variables. However, I expect that very large sample sizes $T$ would be required to apply this method to more than two groups of variables.

We propose the following data analytic technique. 
\begin{enumerate}
\item Given two lists (more than two should also be possible) $X, Y, \ldots$ and $x, y, \ldots$ eliminate all variables not in either list.
\item Construct ``filtered Curto-Itskov complex'' (Ellis and  Klein \cite{spEaK14.ConcurTopolfMRI}).
\item At each ``frequency level'' (``frame''), $M$, project down to ``upper case'' and ``lower case'' complexes $K$ and $L$.
\item  Construct the join $K \ast L$.
\item Compute the persistent homology of the two level filtration $M \hookrightarrow K \ast L$.
\item Any class that has a lifespan of 2 indicates possible independence.
\end{enumerate}

(There is a persistence angle that the preceding does not cover. Hopefully, a future draft of this paper and the software will include that.) Using this method I have found interesting things in real data. (Specifics will be included in a future draft of this paper.)

\section{Simulations}
Applying the method described in section \ref{S:stat.interp} to simulated data I have found that it works! (Specifics will be included in a future draft of this paper.)

\section{K\"unneth theorem with field coefficients} \label{S:Kunneth.with.field.coefs}
Here we prove the following. 

  \begin{theorem}[K\"unneth theorem with field coefficients]  \label{E:Kunneth.with.field}
Let $F$ be a field. Let $\mcl{C}$ be a chain complex $\cdots \rightarrow C_{p+1} \xrightarrow{\partial_{p+1}} C_{p} \rightarrow \cdots$ consisting of finite dimensional vector spaces over $F$ with linear boundary operator $\partial$. Similarly, let $\mcl{D}$ be a  chain complex $\cdots \rightarrow D_{q+1} \xrightarrow{\bar{\partial}_{q+1}} D_{q} \rightarrow \cdots$ consisting of finite dimensional vector spaces over $F$ with linear boundary operator $\bar{\partial}$. Let $\mcl{C} \otimes_{F} \mcl{D}$ be the tensor product of chain complexes. This means the following. The chain complex in dimension $m$ is $\bigoplus_{k+\ell = m} C_{k} \otimes_{F} D_{\ell}$ and the boundary operator is $\underline{\partial} (c_{k} \otimes d_{\ell}) = (\partial c_{k}) \otimes d_{\ell} + (-1)^{k} c_{k} \otimes (\bar{\partial} d_{\ell})$, where $c_{k} \in C_{k}$ and $d_{\ell} \in D_{\ell}$. Let $p \in \ZZ$. Then the map
	\begin{multline}  \label{E:Kunneth.map}
		Kun: \{ w_{k} \} \otimes \{ z_{\ell} \} \mapsto \{ w_{k} \otimes z_{\ell} \} \text{ defines an isomorphism } \\
			\bigoplus_{k+\ell = p} H_{k}(\mcl{C}; F) \otimes H_{\ell}(\mcl{D}; F) \to H_{p} (\mcl{C} \otimes \mcl{D}; F).
	\end{multline}
Here, $w_{k}$ and $z_{\ell}$ are cycles in $C_{k}$, $D_{\ell}$, resp.
  \end{theorem}

  \begin{proof}
The context in which we work are the graded vector spaces $\Gamma := \bigoplus_{k} C_{k}$ and $\Delta := \bigoplus_{\ell} D_{\ell}$ and their bi-graded tensor product. Thus, if $(k, \ell) \neq (k', \ell')$, even if $k+\ell = k'+\ell'$, we have 
$(C_{k} \otimes D_{\ell}) \cap (C_{k'} \otimes D_{\ell'}) = \{ 0 \}$. First, we show that the map defined in \eqref{E:Kunneth.map} is well-defined. Let $W_{k} \subset C_{k}$ be the subset of cycles in $C_{k}$. Let $Z_{\ell}$ be the subspace of $D_{\ell}$ consisting of cycles. Suppose $w \in W_{k}$, $e \in C_{k+1}$, $z \in Z_{\ell}$, and $f \in D_{\ell+1}$. Then
	\begin{equation*}
		\{w + \partial e\} \otimes \{z + \bar{\partial} f \} = (w \otimes z) 
		  + \underline{\partial} \bigl[ (-1)^{k} (w \otimes f) + (e \otimes z) + (e \otimes \bar{\partial} f)  \bigr].
	\end{equation*}

Let $u^{k}_{1}, \ldots, u^{k}_{m_{k}}$ be a basis 
for $C_{k}$ and $v^{\ell}_{1}, \ldots, v^{\ell}_{n_{\ell}}$ be a basis for $D_{\ell}$. Then we can write each $c_{k} \in C_{k}$ as a $m_{k}$-tuple (row vector) in $F^{m_{k}}$ and we can write each $d_{\ell} \in C_{\ell}$ as a $n_{\ell}$-tuple (row vector) in $F^{n_{\ell}}$. In light of this, we also use the symbol $\partial_{k}$ to mean the $m_{k} \times m_{k-1}$ matrix s.t.\ that matrix multiplication $c_{k} \partial_{k} \in C_{k-1}$ gives the boundary of $c_{k}$. 

Let $\zeta_{k} := \dim W_{k}$. Let $i_{k} := m_{k} - \zeta_{k}$. Thus, $i_{k}$ is the dimension of the boundary subspace of $C_{k-1}$. Therefore, $i_{k}  \leq \zeta_{k-1}$. Similarly, let $\phi_{\ell} := \dim Z_{\ell}$. Let $j_{\ell} := n_{\ell} - \phi_{\ell}$. Thus, $j_{\ell}$ is the dimension of the boundary subspace of $D_{\ell-1}$. Therefore, $j_{\ell}  \leq \phi_{\ell-1}$. In summary,
	\begin{equation}  \label{E:ijs.+1.leq.cycle.dim}
		i_{k+1} \leq \zeta_{k} \text{ and } j_{\ell+1} \leq \phi_{\ell}.
	\end{equation} 
Argue as in the proof of Stoll and Wong \cite[Proof of Theorem 2.1, p.\ 99]{rrSetW68.LinearAlgebra} as follows. Let $S \subset C_{k}$ be a complementary subspace of $W_{k}$. I.e., $C_{k} = W_{k} \oplus S$. Thus, $\dim S = i_{k}$. Then the restriction of $\partial_{k} \vert_{S}$ is an isomorphism onto $B_{k-1}$. Let $j : B_{k-1} \to S$ be the inverse of that isomorphism. Hence, $\partial_{k} \circ j$ is the identity on $B_{k-1}$.

Define a chain complex $\{ C_{k}', \partial_{k}' \}$ by $C_{k}' := W_{k} \oplus B_{k-1}$ and $\partial_{k}'(w, b) = (b,0)$ 
($(w,b) \in C_{k}'$). Define $f_{k} : C_{k}' \to C_{k}$ by $f_{k}(w, b) := w + j(b)$. It is easy to see that $f := \{ f_{k} \}$ is a chain map. Notice that $f_{k}$ is an isomorphism for all $k$, so $f$ is a chain isomorphism. Therefore, the homomorphism $f_{\ast}$ in homology is an isomorphism.
We define a complex $\mcl{D}'$ and chain isomorphism $\mcl{D}' \xrightarrow{g} \mcl{D}$ similarly. 

Write $H_{\ast}(\cdot) = H_{\ast}(\cdot; F)$. Similarly, write $\otimes = \otimes_{F}$. It is not hard to see that 
$\bigl( \{ w \}, \{ z \} \bigr) \mapsto \{ w \otimes z \}$ is bilinear map 
$H_{\ast}(\mcl{C}) \times H_{\ast}(\mcl{D}) \to H_{\ast}(\mcl{C} \otimes \mcl{D})$. Therefore, by definition of tensor product 
(Lang \cite[p.\ 408]{sL65.Algebra}) $\eta : \{ w \} \otimes \{ z \} \mapsto \{ w \otimes z \}$ defines a linear homomorphism 
$H_{\ast}(\mcl{C}) \otimes H_{\ast}(\mcl{D}) \to H_{\ast}(\mcl{C} \otimes \mcl{D})$. Define $\eta' : H_{\ast}(\mcl{C}') \otimes H_{\ast}(\mcl{D}') \to H_{\ast}(\mcl{C}' \otimes \mcl{D}')$ similarly. It is not hard to see that the following commutes.
	\begin{equation*} \label{???}
		 \begin{CD}
		   H_{\ast}(\mcl{C}') \otimes H_{\ast}(\mcl{D}') @>{\eta'}>> H_{\ast}(\mcl{C}' \otimes \mcl{D}') \\
		   @V{f_{\ast} \otimes g_{\ast}}V{\isomto}V  @V{\isomto}V{(f \otimes g)_{\ast}}V \\
		   H_{\ast}(\mcl{C}) \otimes H_{\ast}(\mcl{D}) @>{\eta}>> H_{\ast}(\mcl{C} \otimes \mcl{D}) .
 		\end{CD}
	\end{equation*}
We will prove that $\eta'$ is an isomorphism. That will imply that $\eta$ is an isomorphism and complete the proof of the theorem. For now on we assume $\mcl{C} = \mcl{C}'$ and $\mcl{D} = \mcl{D}'$. 

WLOG we may assume that $u^{k}_{1}, \ldots, u^{k}_{\zeta_{k}}$ is a basis for the space of cycles in $C_{k}$, and in fact, we may assume that $u^{k}_{1}, \ldots, u^{k}_{i_{k+1}}$ is a basis for the boundary space in $C_{k}$. Similarly for $D_{\ell}$.
Therefore, 
	\begin{multline}  \label{E:homology.bases}
		u^{k}_{i_{k+1}+1} + Boundaries, \ldots, u^{k}_{\zeta_{k}}  + Boundaries \text{ span } H_{k}(\mcl{C}) \text{ and } \\
		  v^{\ell}_{j_{\ell+1}+1} + Boundaries, \ldots, v^{\ell}_{\phi_{\ell}} + Boundaries \text{ span } H_{\ell}(\mcl{D}).
	\end{multline}
	
W.r.t.\ these bases, the matrices $\partial_{k}$ and $\bar{\partial}_{\ell}$ have the following forms, acting on the right. (We use the convention that superscripts of the form ${}^{\cdot \times \cdot}$ indicate the dimension of matrices.)
	\begin{equation}   \label{E:partial.bar.partial}
		\partial_{k} = 
			\begin{pmatrix}
				0^{\zeta_{k} \times i_{k}} & 0^{\zeta_{k} \times (m_{k-1} - i_{k})} \\
				I_{i_{k}}^{i_{k} \times i_{k}} & 0^{i_{k} \times (m_{k-1} - i_{k})}
			\end{pmatrix} ^{m_{k} \times m_{k-1}} \text{ and } \quad
					\bar{\partial}_{\ell} = 
			\begin{pmatrix}
				0^{\phi_{\ell} \times j_{\ell}} & 0^{\phi_{\ell} \times (n_{\ell-1} - j_{\ell})} \\
				I_{j_{\ell}} & 0^{j_{\ell} \times (n_{\ell-1} - j_{\ell})}
			\end{pmatrix} ^{n_{\ell} \times n_{\ell-1}} ,
	\end{equation}	 
where $I_{i_{k}}$ the ${i_{k} \times i_{k}}$ identity matrix.
	
$u^{k}_{i} \otimes v^{\ell}_{j}$ ($i = 1, \ldots, m_{k}$; $j = 1, \ldots, n_{\ell}$) is a basis of $C_{k} \otimes D_{\ell}$. Let 
$c_{k} = x_{1} u^{k}_{1} + \cdots + x_{m_{k}} u^{k}_{m_{k}} \in C_{k}$. 
Let $d_{\ell} = y_{1} v^{\ell}_{1} + \cdots + y_{m_{\ell}} v^{\ell}_{n_{\ell}} \in D_{\ell}$.
(Here $x_{1}, \ldots, x_{m_{k}}, y_{1}, \ldots, y_{n_{\ell}} \in F$, of course.) 
Hence, $c_{k} \otimes d_{\ell} = \sum_{i=1}^{m_{k}} \sum_{j=1}^{n_{\ell}} x_{i} y_{j} \, u^{k}_{i} \otimes v^{\ell}_{j}$. Thus, the coefficient of $u^{k}_{i} \otimes v^{\ell}_{j}$ in $c_{k} \otimes d_{\ell}$ is just the $(i,j)^{th}$ entry in the matrix 
$c_{k}^{T} d_{\ell} = ( x_{i} y_{j} )^{m_{k} \times n_{\ell}}$.    

Thus, a typical element, $e_{k, \ell}$ of $C_{k} \otimes D_{\ell}$ is represented as an $m_{k} \times n_{\ell}$ matrix, $E_{k, \ell}$. Conversely, \emph{an arbitrary $m_{k} \times n_{\ell}$ matrix (with entries in $F$) represents an element of $C_{k} \otimes D_{\ell}$.} (\emph{Pf:} Let $E = (e_{ij})^{m_{k} \times n_{\ell}}$. Then $E$ represents 
$\sum_{i=1}^{m_{k}} \sum_{j=1}^{n_{\ell}} e_{ij} \, u^{k}_{i} \otimes v^{\ell}_{j}$.) In this way we get an isomorphism from $C_{k} \otimes D_{\ell}$ to the space, $M_{k,\ell}$, of $m_{k} \times n_{\ell}$ matrices with entries in $F$.  From now on we identify $C_{k}$ with $F^{m_{k}}$, $D_{\ell}$ with $F^{n_{\ell}}$ and $C_{k} \otimes D_{\ell}$ with $M_{k,\ell}$. Thus, a typical element of $\mcl{C} \otimes \mcl{D}$ can be written
	\begin{equation*}
		\sum_{k} \sum_{\ell} E_{k, \ell}^{m_{k} \times n_{\ell}},
	\end{equation*}
where almost all $E_{k, \ell}$'s are 0. We have
	\begin{equation*}
		\underline{\partial} E_{k, \ell} = \partial^{T} E_{k, \ell} + (-1)^{k} E_{k, \ell} \bar{\partial}. 
	\end{equation*}
It is convenient to write 
	\begin{equation}  \label{E:general.form.of.tensor}
		E_{k,\ell} = 
			\begin{pmatrix}
				(E_{k,\ell}^{1,1})^{\zeta_{k} \times \phi_{\ell}} & (E_{k,\ell}^{1,2})^{\zeta_{k} \times j_{\ell}} \\
				(E_{k,\ell}^{2,1})^{i_{k} \times \phi_{\ell}} & (E_{k,\ell}^{2,2})^{i_{k} \times j_{\ell}}
			\end{pmatrix} .
	\end{equation}
Then
	\begin{multline} \label{E:general.form.of.boundary}
		\underline{\partial} E_{k,\ell} = 
			\begin{pmatrix}
			   (E_{k,\ell}^{2,1})^{i_{k} \times \phi_{\ell}} & (E_{k,\ell}^{2,2})^{i_{k} \times j_{\ell}}  \\
			   0^{(m_{k-1}-i_{k}) \times \phi_{\ell}} & 0^{(m_{k-1}-i_{k}) \times j_{\ell}}
			\end{pmatrix}^{m_{k-1} \times n_{\ell}} \\
			  \oplus (-1)^{k} 
			\begin{pmatrix}
			   (E_{k,\ell}^{1,2})^{\zeta_{k} \times j_{\ell}} & 0^{\zeta_{k} \times (n_{\ell-1} - j_{\ell})} \\
			   (E_{k,\ell}^{2,2})^{i_{k} \times j_{\ell}} & 0^{i_{k} \times (n_{\ell-1} - j_{\ell})}
			\end{pmatrix}^{m_{k} \times n_{\ell-1}} \\
			 \in (C_{k-1} \otimes D_{\ell}) \oplus (C_{k} \otimes D_{\ell-1}) \subset (\mcl{C} \otimes \mcl{D})_{k+\ell-1},
	\end{multline}
where we separate the matrices by ``$\oplus$'' as a reminder that they have different bi-degrees. 

\subsection{Cycles}  \label{SS:cycles}
First, we prove that any cycle in $\mcl{C} \otimes \mcl{D}$ is homologous to an element of $\mcl{W} \otimes \mcl{Z}$, where 
$\mcl{W} = \bigoplus_{k} W_{k}$ and $\mcl{Z} = \bigoplus_{\ell} Z_{\ell}$. First of all note that $c_{k} = (x_{1}, \ldots, x_{m_{k}}) \in C_{k}$ is in $W_{k}$ if and only if $x_{\zeta_{k}+1} = \cdots = x_{m_{k}} = 0$. Similarly, $d_{\ell} = (y_{1}, \ldots, y_{n_{\ell}}) \in D_{\ell}$ is in $Z_{\ell}$ if and only if $y_{\phi_{\ell}+1} = \cdots = y_{n_{\ell}} = 0$. Therefore, 
an element of the $(k, \ell)^{th}$ component of $\mcl{W} \otimes \mcl{Z}$ can be written in the form \eqref{E:general.form.of.tensor} with $E_{k,\ell}^{1,2}$, $E_{k,\ell}^{2,1}$, and $E_{k,\ell}^{2,2}$ all 0. Conversely, by \eqref{E:general.form.of.boundary}, any such element is a cycle of $\mcl{C} \otimes \mcl{D}$. We show, then, that any cycle in $\mcl{C} \otimes \mcl{D}$ is at least homologous to cycle of that form.

Consider a matrix $E_{k+1,\ell+1} \in C_{k+1} \otimes D_{\ell+1}$ with $E_{k+1,\ell+1}^{1,2}$, $E_{k+1,\ell+1}^{2,1}$, and $E_{k+1,\ell+1}^{1,1}$ all 0. Then, by \eqref{E:general.form.of.boundary},
	\begin{multline} \label{E:E22.only.boundary}
		\underline{\partial} E_{k+1,\ell+1} = 
			\begin{pmatrix}
			   0^{i_{k+1} \times \phi_{\ell+1}} & (E_{k+1,\ell+1}^{2,2})^{i_{k+1} \times j_{\ell+1}}  \\
			   0^{(m_{k}-i_{k+1}) \times \phi_{\ell+1}} & 0^{(m_{k}-i_{k+1}) \times j_{\ell+1}}
			\end{pmatrix}^{m_{k} \times n_{\ell+1}} \\
			  + (-1)^{k+1} 
			\begin{pmatrix}
			   0^{\zeta_{k+1} \times j_{\ell+1}} & 0^{\zeta_{k+1} \times (n_{\ell} - j_{\ell+1})} \\
			   (E_{k+1,\ell+1}^{2,2})^{i_{k+1} \times j_{\ell+1}} & 0^{i_{k+1} \times (n_{\ell} - j_{\ell+1})}
			\end{pmatrix}^{m_{k+1} \times n_{\ell}} \\
			 \in (C_{k} \otimes D_{\ell+1}) \oplus (C_{k+1} \otimes D_{\ell}).
	\end{multline}

Next, let $z = \sum_{k+\ell=p} A_{k,\ell}$ (finite sum) be a cycle in $(\mcl{C} \otimes \mcl{D})_{p}$. 
Then $0 = \underline{\partial} z = \sum_{k,\ell} \bigr[ \partial^{T} A_{k,\ell} + (-1)^{k} A_{k,\ell} \bar{\partial} \bigl]$. For every $k$ and $\ell$, with $k + \ell = p$ the $(k,\ell)$-component of $\underline{\partial} z$ must be 0. I.e., for every $k$ and $\ell$ s.t.\ $k + \ell = p$, we must have $\partial^{T} A_{k+1,\ell} + (-1)^{k} A_{k,\ell+1} \bar{\partial} = 0$. Thus, by \eqref{E:general.form.of.boundary},
	\begin{multline} \label{E:pieces.of.boundary}
		\begin{pmatrix}
			   (A_{k+1,\ell}^{2,1})^{i_{k+1} \times \phi_{\ell}} & (A_{k+1,\ell}^{2,2})^{i_{k+1} \times j_{\ell}}  \\
			   0^{(m_{k}-i_{k+1}) \times \phi_{\ell}} & 0^{(m_{k}-i_{k+1}) \times j_{\ell}} \\
		\end{pmatrix}^{m_{k} \times n_{\ell}} \\
		  + (-1)^{k} 
		\begin{pmatrix}
		   (A_{k,\ell+1}^{1,2})^{\zeta_{k} \times j_{\ell+1}} & 0^{\zeta_{k} \times (n_{\ell} - j_{\ell+1})} \\
		   (A_{k,\ell+1}^{2,2})^{i_{k} \times j_{\ell+1}} & 0^{i_{k} \times (n_{\ell} - j_{\ell+1})}
		\end{pmatrix}^{m_{k} \times n_{\ell}} 
			   = 0.
	\end{multline}
We immediately conclude that $A_{k+1,\ell}^{2,2}$ and $A_{k,\ell+1}^{2,2}$ are both 0.

By \eqref{E:ijs.+1.leq.cycle.dim} we may write
	\begin{equation*}
	  A_{k+1,\ell}^{2,1} = \bigl( G_{k+1,\ell}^{i_{k+1} \times j_{\ell+1}}, H_{k+1,\ell}^{i_{k+1} \times (\phi_{\ell} - j_{\ell+1}} ) \text{ and }
	  A_{k,\ell+1}^{1,2} = 
	    \begin{pmatrix}
		   G_{k,\ell+1}^{i_{k+1} \times j_{\ell+1}} \\
		   H_{k,\ell+1}^{(\zeta_{k} - i_{k+1}) \times j_{\ell+1}}
	    \end{pmatrix} .
	\end{equation*}
From \eqref{E:pieces.of.boundary}, we see that $H_{k+1,\ell}$ and $H_{k,\ell+1}$ are both 0 and 
$G_{k,\ell+1}  = (-1)^{k+1} G_{k+1,\ell}$.

Thus,
	\begin{multline*}  
		A_{k+1,\ell} = 
			\begin{pmatrix}
				(A_{k+1,\ell}^{1,1})^{\zeta_{k+1} \times \phi_{\ell}} & (A_{k+1,\ell}^{1,2})^{\zeta_{k+1} \times j_{\ell}} \\
				( G_{k+1,\ell}, 0^{i_{k+1} \times (\phi_{\ell} - j_{\ell+1})} ) 
				     & 0^{i_{k+1} \times j_{\ell}}
			\end{pmatrix} \text{ and } \\
		A_{k,\ell+1} = 
			\begin{pmatrix}
				(A_{k,\ell+1}^{1,1})^{\zeta_{k} \times \phi_{\ell+1}} & 
				     \begin{pmatrix}
					   (-1)^{k+1} G_{k+1,\ell} \\
					   0^{(\zeta_{k} - i_{k+1}) \times j_{\ell+1}}
				    \end{pmatrix}  \\
				(A_{k,\ell+1}^{2,1})^{i_{k} \times \phi_{\ell+1}} & 0^{i_{k} \times j_{\ell+1}}
			\end{pmatrix}
	\end{multline*}

Now let $E_{k+1,\ell+1} \in C_{k+1} \otimes D_{\ell+1}$ be as in \eqref{E:E22.only.boundary} with $E_{k+1,\ell+1}^{2,2} = G_{k+1,\ell}$. Let $E_{k+1,\ell+1}^{2,2} = G_{k+1,\ell}$. Let $\tilde{A}_{k+1,\ell}$ be just like $A_{k+1,\ell}$, but with $\tilde{A}_{k+1,\ell}^{2,1} = 0$. Let $\tilde{A}_{k,\ell+1}$ be just like $A_{k,\ell+1}$, but with $\tilde{A}_{k,\ell+1}^{1,2} = 0$. Then 
	\begin{equation*}
		A_{k+1,\ell} + A_{k,\ell+1} = \tilde{A}_{k+1,\ell} + \tilde{A}_{k,\ell+1} + \underline{\partial} E_{k+1,\ell+1}.
	\end{equation*}

Thus, we may assume, up to homology, 
	\begin{equation*}  
		A_{k+1,\ell} = 
			\begin{pmatrix}
				(A_{k+1,\ell}^{1,1})^{\zeta_{k+1} \times \phi_{\ell}} & (A_{k+1,\ell}^{1,2})^{\zeta_{k+1} \times j_{\ell}} \\
				0^{i_{k+1} \times \phi_{\ell}} & 0^{i_{k+1} \times j_{\ell}}
			\end{pmatrix} \text{ and } 
		A_{k,\ell+1} = 
			\begin{pmatrix}
				(A_{k,\ell+1}^{1,1})^{\zeta_{k} \times \phi_{\ell+1}} & 0^{\zeta_{k} \times j_{\ell+1}} \\
				(A_{k,\ell+1}^{2,1})^{i_{k} \times \phi_{\ell+1}} & 0^{i_{k} \times j_{\ell+1}}
			\end{pmatrix}
	\end{equation*}
Playing the same games with $\partial^{T} A_{k+2,\ell-1} + (-1)^{k+1} A_{k+1,\ell} \bar{\partial} = 0$ and $\partial^{T} A_{k,\ell+1} + (-1)^{k-1} A_{k-1,\ell+2} \bar{\partial} = 0$ we see that we may take $A_{k+1,\ell}^{1,2} = 0$ and $A_{k,\ell+1}^{2,1} = 0$.

Hence, finally, up to homology, 
	\begin{equation*}  
		A_{k+1,\ell} = 
			\begin{pmatrix}
				(A_{k+1,\ell}^{1,1})^{\zeta_{k+1} \times \phi_{\ell}} & 0^{\zeta_{k+1} \times j_{\ell}} \\
				0^{i_{k+1} \times \phi_{\ell}} & 0^{i_{k+1} \times j_{\ell}}
			\end{pmatrix} \text{ and } 
		A_{k,\ell+1} = 
			\begin{pmatrix}
				(A_{k,\ell+1}^{1,1})^{\zeta_{k} \times \phi_{\ell+1}} & 0^{\zeta_{k} \times j_{\ell+1}} \\
				0^{i_{k} \times \phi_{\ell+1}} & 0^{i_{k} \times j_{\ell+1}}
			\end{pmatrix}
	\end{equation*}
But as noted above these matrices represent elements of $\mcl{W} \otimes \mcl{Z}$, as desired.
	
\subsection{Basis of $H_{k}(\mcl{C}) \otimes H_{\ell}(\mcl{D})$.}	
Let $k, \ell, p \in \ZZ$ with $k + \ell = p$. Let $A_{k, \ell}$ be a cycle of bi-degree $(k, \ell)$. From subsection \ref{SS:cycles}, we know that we may assume 
	\begin{equation*}
		A_{k,\ell} = 
			\begin{pmatrix}
				(A_{k,\ell}^{1,1})^{\zeta_{k} \times \phi_{\ell}} & 0^{\zeta_{k} \times j_{\ell}} \\
				0^{i_{k} \times \phi_{\ell}} & 0^{i_{k} \times j_{\ell}}
			\end{pmatrix}^{m_{k} \times n_{\ell}}.
	\end{equation*}
Write
	\begin{equation*}
		A_{k,\ell}^{1,1} = 
			\begin{pmatrix}
				\alpha^{i_{k+1} \times j_{\ell+1}} & \beta^{i_{k+1} \times (\phi_{\ell} - j_{\ell+1})} \\
				\gamma^{(\zeta_{k} - i_{k+1}) \times j_{\ell+1}} & \delta^{(\zeta_{k} - i_{k+1}) \times (\phi_{\ell} - j_{\ell+1})}
			\end{pmatrix}.
	\end{equation*}
(See \eqref{E:ijs.+1.leq.cycle.dim}.)

For any $E_{k+1,\ell}$ and $E_{k,\ell+1}$ of indicated bi-degree, both $\underline{\partial} E_{k+1,\ell}$ and $\underline{\partial} E_{k,\ell+1}$ have components of bi-degree $(k, \ell)$ and only such $E$'s can do so. But we see from \eqref{E:general.form.of.boundary}, that neither of those components have non-zero entries in the positions covered by $\delta$. Hence, the entries in $\delta$ correspond to basis elements of $H_{k+\ell}( \mcl{C} \otimes \mcl{D} )$. 
	
Conversely, let
	\begin{equation*}
		\tilde{A}_{k,\ell}^{1,1} = 
			\begin{pmatrix}
				\alpha^{i_{k+1} \times j_{\ell+1}} & \beta^{i_{k+1} \times (\phi_{\ell} - j_{\ell+1})} \\
				\gamma^{(\zeta_{k} - i_{k+1}) \times j_{\ell+1}} & 0^{(\zeta_{k} - i_{k+1}) \times (\phi_{\ell} - j_{\ell+1})}
			\end{pmatrix},
	\end{equation*}
let
	\begin{equation*}  
		E_{k+1,\ell} = 
			\begin{pmatrix}
				0^{\zeta_{k+1} \times \phi_{\ell}} & 0^{\zeta_{k+1} \times j_{\ell}} \\
				(E_{k+1,\ell}^{2,1})^{i_{k+1} \times \phi_{\ell}} & 0^{i_{k+1} \times j_{\ell}}
			\end{pmatrix}, 
	       \text{ where, } E_{k+1,\ell}^{2,1} = 
			\begin{pmatrix}
				\alpha^{i_{k+1} \times j_{\ell+1}} & \beta^{i_{k+1} \times (\phi_{\ell} - j_{\ell+1})} 
			\end{pmatrix},
	\end{equation*}
and let
	\begin{equation*}  
		E_{k,\ell+1} = 
			\begin{pmatrix}
				0^{\zeta_{k} \times \phi_{\ell+1}} & (E_{k,\ell+1}^{1,2})^{\zeta_{k} \times j_{\ell+1}} \\
				0^{i_{k} \times \phi_{\ell+1}} & 0^{i_{k} \times j_{\ell+1}}
			\end{pmatrix}, 
	       \text{ where } E_{k,\ell+1}^{1,2} = (-1)^{k}
			\begin{pmatrix}
				0^{i_{k+1} \times \phi_{j+1}} \\ 
				\gamma^{(\zeta_{k} - i_{k+1}) \times j_{\ell+1}} 
			\end{pmatrix}.
	\end{equation*}
Then, by \eqref{E:general.form.of.boundary}, 
	\begin{equation*} 
		\underline{\partial} E_{k+1,\ell} +  \underline{\partial} E_{k,\ell+1} = 
			\begin{pmatrix}
				(\tilde{A}_{k,\ell}^{1,1})^{\zeta_{k} \times \phi_{\ell}} & 0^{\zeta_{k} \times j_{\ell}} \\
				0^{i_{k} \times \phi_{\ell}} & 0^{i_{k} \times j_{\ell}}
			\end{pmatrix}^{m_{k} \times n_{\ell}}.
	\end{equation*}
This shows that $A_{k,\ell}$ is homologous to
	\begin{equation*}
	   F_{k, \ell} :=
			\begin{pmatrix}
				(G_{k,\ell}^{1,1})^{\zeta_{k} \times \phi_{\ell}} & 0^{\zeta_{k} \times j_{\ell}} \\
				0^{i_{k} \times \phi_{\ell}} & 0^{i_{k} \times j_{\ell}}
			\end{pmatrix}^{m_{k} \times n_{\ell}}
		\text{ with } G_{k,\ell}^{1,1} = 
			\begin{pmatrix}
				0^{i_{k+1} \times j_{\ell+1}} & 0^{i_{k+1} \times (\phi_{\ell} - j_{\ell+1})} \\
				0^{(\zeta_{k} - i_{k+1}) \times j_{\ell+1}} & \delta^{(\zeta_{k} - i_{k+1}) \times (\phi_{\ell} - j_{\ell+1})}
			\end{pmatrix}.
	\end{equation*}
 
For $a = 1, ..., m_{k}$, $b = 1, \ldots, n_{\ell}$, let $1_{k, \ell}^{a,b}$ be the $m_{k} \times n_{\ell}$ matrix all of whose entries are 0 except at $(a,b)$ where there is a 1. We have proved that a basis for $H_{p}( \mcl{C} \otimes \mcl{D} )$ is
	\begin{equation} \label{E:basis.of.Hp}
		\text{All } \{1_{k, \ell}^{a,b}\} \text{ with } k + \ell = p; \; a = i_{k}+1, \ldots, \zeta_{k}; 
		  \text{ and } b = j_{\ell}+1, \ldots, \phi_{\ell} .
	\end{equation}

 By \eqref{E:homology.bases}, a basis for $H_{k}(\mcl{C}) \otimes H_{\ell}(\mcl{D})$ is 
	 \begin{equation}  \label{E:basis.of.Hk.ox.Hell}
		\{ u^{k}_{a} \} \otimes \{ v^{\ell}_{b} \} \quad
		  (a = i_{k+1}+1, \ldots, \zeta_{k}, b = j_{\ell+1}+1, \ldots, \phi_{\ell}). 
	\end{equation}
By convention, $u^{k}_{a}$ is the $m_{k}$ row vector that is all 0 except for a 1 in position $a$. Similarly for $v^{\ell}_{b}$. Recall the definition of the function $Kun$ in the statement of the theorem. We thus have 
	\begin{equation*}
		Kun \bigl( \{ u^{k}_{a} \} \otimes \{ v^{\ell}_{b} \} \bigr) = (u^{k}_{a})^{T}v^{\ell}_{b} = \{1_{k, \ell}^{a,b}\}. 
	\end{equation*}
Letting $a$ and $b$ vary as in \eqref{E:basis.of.Hk.ox.Hell} and letting $k, \ell$ vary over $k + \ell = p$ and comparing to \eqref{E:basis.of.Hp}, we see that $Kun$ is an isomorphism.
  \end{proof}
  
\section{Acknowlegements}
Steve Ferry helped with some technicalities.


\providecommand{\bysame}{\leavevmode\hbox to3em{\hrulefill}\thinspace}
\providecommand{\MR}{\relax\ifhmode\unskip\space\fi MR }
\providecommand{\MRhref}[2]{%
  \href{http://www.ams.org/mathscinet-getitem?mr=#1}{#2}
}
\providecommand{\href}[2]{#2}

\end{document}